\documentclass[a4paper,11pt]{amsart}
\usepackage{amssymb}
\usepackage[english]{babel}
\usepackage[latin1]{inputenc}
\usepackage{graphicx}

\theoremstyle{plain}
\newcommand{\ep}{\epsilon}
\newcommand{\R}{\mathbb R}

\begin{document}
\title[A boundary point lemma]{A boundary point lemma for Black-Scholes type operators}
\author[Erik Ekström and Johan Tysk]{Erik Ekström$^1$ and Johan Tysk$^{2,3}$}
\email{ekstrom@maths.manchester.ac.uk, Johan.Tysk@math.uu.se}
\subjclass[2000]{Primary 35K65, 35B50; Secondary 91B28}
\keywords{Hopf boundary point lemma, parabolic equations, degenerate
equations, Black-Scholes equation}
\thanks{$^1$ School of Mathematics, The University of Manchester, Sackville Street,
Manchester M60 1QD, UK}
\thanks{$^2$ Department of Mathematics, Uppsala University, Box 480, SE-75106 Uppsala, Sweden}
\thanks{$^3$ Partially supported by the Swedish Research Council (VR)}

\newtheorem{theorem}{Theorem}[section]
\newtheorem{lemma}[theorem]{Lemma}
\newtheorem{corollary}[theorem]{Corollary}
\newtheorem{proposition}[theorem]{Proposition}
\newtheorem{definition}[theorem]{Definition}
\newtheorem{hypothesis}[theorem]{Hypothesis}

\newenvironment{example}[1][Example]{\begin{trivlist}
\item[\hskip \labelsep {\bf Example}]}{\end{trivlist}}
\newenvironment{remark}[1][Remark]{\begin{trivlist}
\item[\hskip \labelsep {\bf Remark}]}{\end{trivlist}}

\begin{abstract}
We prove a sharp version of the Hopf boundary point lemma for Black-Scholes type equations.
We also investigate the existence and the regularity of the spatial derivative of the solutions
at the spatial boundary.
\end{abstract}

\maketitle

\section{Introduction}

In this paper we study a parabolic initial-boundary value problem motivated by
applications involving modeling of non-negative quantities described by diffusion processes
that are absorbed at 0. To describe the problem in a financial setting, consider a market
consisting of a risk-free asset
\[dB=r(t)B\,dt\]
where $r$ is a deterministic function,
and $n$ risky assets with non-negative prices $X_i$, $i=1,...,n$, modeled by diffusion processes.
According to standard arbitrage theory, to price options maturing at some future time $T$
written on such risky assets,
we need not specify the expected rate of return of the assets, compare for example Theorem~7.8 and
Proposition~7.9 in \cite{B} and Proposition~2.2.3 in \cite{KS}. Instead
one specifies the price dynamics under a so-called risk-neutral measure, under which the
expected rate of return of all assets equals the interest rate $r(t)$.
For the sake of notational convenience, however,
replacing $X_i(t)$ by $X_i(t)\exp\{\int_t^T r(s)\,ds\}$, we obtain processes with drift 0.
Financially this corresponds to quoting the assets in terms of bonds maturing at time $T$.
By abuse of notation, we denote also these new processes by $X_i$.
Thus we model the price dynamics of the $i$th asset by
\[dX_i=\sum_{j=1}^n\alpha_{ij}(X_i,t)\,dW_j\hspace{10mm}X_i(t)=x_i.\]
Here $W_j$, $j=1,...,n$, are independent standard Brownian motions, and
$\alpha=(\alpha_{ij})^n_{i,j=1}$ is for each $x$ and $t$ an $n\times n$-matrix.
We assume that the components of $X=(X_1,...,X_n)$ are absorbed at 0, i.e. if some component at some
point reaches 0, then it remains 0 forever. Therefore we let $\alpha_{ij}(0,t)=0$ for all
$i$ and $j$. We also assume that the rank of $\alpha$ is equal to the number of non-zero
spatial coordinates (this assumption is connected to the completeness of the model; more precisely, under
this assumption every option has a unique arbitrage-free price, compare Corollary~12.2.6 in \cite{O} or
Theorem~1.6.6 in \cite{KS}).
Further assumptions on $\alpha$ are specified below.
Given a continuous contract function $g:[0,\infty)^n\to\R$ of at most polynomial growth,
the value $U(x,t)$ at time $t$ of an option paying $g(X(T))$ at time $T$ is
\[U(x,t):=E_{x,t}g\big(X(T)\big),\]
where the indices indicate that $X(t)=x$. Alternatively, the function $u(x,t):=U(x,T-t)$ solves
the Black-Scholes equation
\begin{equation}
\label{u}
\left\{ \begin{array}{l}
\mathcal L u=0 \mbox{ for $(x,t)\in(0,\infty)^n\times(0,\infty)$}\\
u(x,0)=g(x),
\end{array}\right.
\end{equation}
where
\[\mathcal Lu=\sum_{i,j=1}^n a_{ij}u_{x_ix_j}-u_t,\]
$a_{ij}=a_{ij}(x,t)$ are the elements of the matrix
\[a(x,t):=\alpha(x,T-t)\alpha^*(x,T-t)/2\]
and
\[u_t:=\frac{\partial u}{\partial t},\hspace{15mm}
u_{x_ix_j}:=\frac{\partial^2u}{\partial x_j\partial x_j}.\]
In order to have uniqueness of solutions to the problem (\ref{u}) in the class of functions
with at most polynomial growth,
one needs (in general) to impose boundary conditions at the faces $\{x_i=0\}$, $i=1,...,n$.
We assume that these boundary conditions are defined inductively by solving the partial
differential equation in lower dimensional faces
(the assumptions on $\alpha$ guarantee that the operator is parabolic in these faces).
One thus starts with solving an ordinary differential equation ($u_t=0$) along
the $t$-axis, then one solves $n$ parabolic equations in the faces spanned by $t$ and one
of the variables $x_1,...,x_n$, and so on. According to Theorems~4.1 and 5.5 in \cite{JT}, this
procedure results in a unique classical solution of polynomial growth to equation (\ref{u}).

In this paper we are concerned with the behavior of the first order derivatives $u_{x_i}$ at the boundary.
We should note that the derivatives $u_{x_i}$ are of special importance in financial applications.
They are the so-called ``deltas'' of the option, and they represent the number of
stocks $X_i$ a hedger should have in a hedging portfolio, compare for example Theorem~8.5 in \cite{B}.

The following theorem, the Hopf boundary point lemma (adapted to our current setting),
is well-known in the theory of parabolic equations, compare for example
\cite{F} (Theorem 2.14, p. 49) or \cite{L} (Lemma 2.6, p. 10).

\begin{theorem}
\label{Hopf}
{\bf (Hopf boundary point lemma.)}
Assume that the differential operator is uniformly parabolic, i.e. that there exists
$\gamma>0$ such that $\xi a(x,t)\xi^*\geq \gamma \vert\xi\vert^2$ for all
$\xi\in\R^n$, $x$ and $t$. Let a point $P^\prime=(0,x^\prime_2,...,x^\prime_n,t_0)$
with $x^\prime_2,...,x^\prime_n,t_0>0$ be given. Assume that in a neighborhood of
$P^\prime$ we have $u(P)>u(P^\prime)$ for all points $P$
with $x_1>0$ and $t\leq t_0$. Then $u_{x_1}(P^\prime)>0$ in the sense that
\[\liminf_{\ep\to 0}\frac{1}{\ep}\Big(u(P^\prime+\ep e_1)-u(P^\prime)\Big)>0,\]
where $e_1$ is the unit vector in the $x_1$-direction.
\end{theorem}

The following well-known example shows that the above theorem is not valid without the
assumption of uniform parabolicity.

\begin{example}
{\bf (The Black-Scholes price of a call option.)}
Let $n=1$, and let $\alpha(x,t)=\sigma x$ and $g(x)=(x-K)^+$ where $\sigma$ and $K$ are positive
constants (here and in the sequel we drop the subscripts of $\alpha_{11}$, $a_{11}$ and $x_1$ if
$n=1$). Then
\[u(x,t)=x\Phi\Big(\frac{\ln (x/K)+\sigma^2t/2}{\sqrt{\sigma^2t}}\Big)-
K\Phi\Big(\frac{\ln (x/K)-\sigma^2t/2}{\sqrt{\sigma^2t}}\Big)\]
where
\[\Phi(z)=\frac{1}{\sqrt{2\pi}}\int_{-\infty}^z\exp\{-u^2/2\}\,du.\]
It is now straightforward to check that the delta of the call option is given by
\[u_x(x,t)=\Phi\Big(\frac{\ln (x/K)+\sigma^2t/2}{\sqrt{\sigma^2t}}\Big),\]
so $u_x(0,t)=0$.
\end{example}

The outline of the present article is as follows.
In Section~\ref{hopf} we provide a Hopf lemma for equations that are not uniformly
parabolic, compare Theorems~\ref{degenhopf} and \ref{degenhopf2}.
To prove these results we need to assume that $a_{11}\geq Cx_1^\beta$ for some $0<\beta<2$. In
view of the above example (in which $a_{11}=\sigma^2x_1^2/2$),
our result can be regarded as a sharp version of the Hopf lemma.
Along the lines of the above example we also show that if $n=1$ and $a\leq Cx^2$, then the result
of the Hopf lemma always fails, compare Theorem~\ref{neghopf}.

In Section~\ref{regularity} we perform further investigations
of the spatial derivative at the boundary. We use preservation of convexity in dimension one
to show that $u_x(0,t)$ always exists finitely for a wide class of initial conditions,
and we provide an example in which $u_x(0,t)$ is discontinuous.

\section{Hopf boundary point lemma}
\label{hopf}
In this section we prove a version of the Hopf boundary point lemma for non-uniformly
parabolic operators.

\begin{hypothesis}
\label{hyp}
The $n\times n$-matrix $\alpha=\big(\alpha_{ij}(x,t)\big)_{i,j=1}^n$ is defined on
$[0,\infty)^n\times (-\infty,T]$, and
\begin{itemize}
\item[(i)]
is continuous in the time variable, and $\alpha$ is also Lipschitz continuous in the spatial variables
on every compact subset of $(0,\infty)^n\times (-\infty,T]$; we also assume the standard growth condition
$\vert\alpha(x,t)\vert\leq C(1+\vert x\vert)$ for some constant $C$;
\item[(ii)]
for each pair $(i,j)$, $\alpha_{ij}$ is a function merely of $x_i$ and $t$ (this is automatic if $n=1$),
and for $x_i=0$ we have $\alpha_{ij}=0$;
\item[(iii)]
the rank of $\alpha$ is equal to the number of non-zero spatial coordinates.
\end{itemize}
If $n=1$ the
local Lipschitz condition may be replaced with a local H\"older(1/2) condition.
\end{hypothesis}

Note that we allow discontinuities of the coefficients at the spatial boundary.
Also note that condition (ii) does of course not rule
out the possibility of dependence between the different assets. It merely means that the instantaneous
covariance $a_{ij}$ between two assets $X_i$ and $X_j$ only depends on time and the present
values of $X_i$ and $X_j$. The condition (iii) is connected to the issue of completeness of the
model, compare the discussion in the introduction.

\begin{theorem}
\label{degenhopf}
Assume Hypothesis~\ref{hyp}, that $u$ satisfies the Black-Scholes equation (\ref{u}),
and that $a_{11}(x_1,t)\geq Cx_1^\beta$ for some constants $C>0$ and $\beta\in[0,2)$.
Let a point $P^\prime=(0,x^\prime_2,...,x^\prime_n,t_0)$
with $t_0>0$ and $x^\prime_i>0$ for $i=2,...,n$ be given.
Assume that in some neighborhood of $P^\prime$ we have that
\[u(x_1,x_2,...,x_n,t)>u(0,x_2,...,x_n,t)\]
for all points with $x_1>0$ and
$t\leq t_0$. Then $u_{x_1}(P^\prime)> 0 $ (in the same sense as described in Theorem~\ref{Hopf}).
\end{theorem}

\begin{proof}
Introduce the function
\[v(x,t):=x_1+x_1^{1+\ep}-\vert t-t_0\vert^N-\sum_{i=2}^n\vert x_i-x^\prime_i\vert^{2N}\]
for some constants $\ep>0$ and $N\geq 1$, both to be chosen later. Let
\[D:=\{ P:v(P)\geq 0,0\leq x_1\leq \eta, t\leq t_0\}\]
for some small constant $\eta>0$.
Then there exists a positive constant $C_1$ (depending on $\eta$) such that
\[-C_1x_1^{1/N}\leq t-t_0\leq 0\]
and
\[-C_1x_1^{1/(2N)}\leq -\vert x_i-x^\prime_i\vert\leq 0,\]
$i=2,...n$, for points in $D$. Consequently,
\begin{eqnarray*}
\mathcal L v &=& \sum_{i,j=1}^n a_{ij}v_{x_ix_j}-v_t
= \sum_{i=1}^n a_{ii}v_{x_ix_i}-v_t\\
&=& a_{11}(1+\ep)\ep x_1^{\ep-1}-N\vert t-t_0\vert^{N-1}-
(4N^2-2N)\sum_{i=2}^na_{ii}\vert x_i-x_i^\prime\vert^{2N-2}\\
&\geq& C(1+\ep)\ep x_1^{\beta+\ep-1}-NC_2x_1^{(N-1)/N} - (4N^2-2N) C_2x_1^{(N-1)/N}
\end{eqnarray*}
in $D$ for some constant $C_2$ satisfying
\[C_2>\max\big\{C_1^{N-1},C_1^{N-2}(n-1)\max_{2\leq i\leq n}\{a_{ii}(P^\prime)\}\big\}.\]
Thus, choosing $\ep$ small and $N$ large so that $\beta+\ep-1<(N-1)/N<1$, it is
clear that $\mathcal L v\geq 0$ in $D$ (at least if $\eta$ is small enough; note that
$C_1$, and thus also $C_2$, can be held fixed when decreasing $\eta$).
The parabolic boundary $\partial_p D$ of $D$ can be written as
$\partial_p D=S_1\cup S_2$ where
\[S_1=\{P: v(P)=0, t\leq t_0, x_1\leq \eta\}\]
and
\[S_2=\{P:v(P)\geq 0,t\leq t_0, x_1=\eta\}.\]
Since, by assumption,
\[u(x_1,x_2,...,x_n,t)>u(0,x_2,...,x_n,t)\mbox{ for }(x_1,x_2,...,x_n,t)\in D\setminus \{P^\prime\}\]
there exists $\delta>0$ such that $u(x_1,x_2,...,x_n,t)-u(0,x_2,...,x_n,t)\geq \delta$ for
$(x_1,x_2,...,x_n,t)\in S_2$.
Since $v$ is bounded by 1 in $D$ (at least if $\eta$ is small enough) it follows that
\begin{equation}
\label{eq1}
u(x_1,x_2,...,x_n,t)-u(0,x_2,...,x_n,t)-\delta v(x_1,x_2,...,x_n,t)\geq 0 \mbox{ on $\partial_p D$}.
\end{equation}
Moreover, applying the differential operator
\[\mathcal L=\sum_{i,j=1}^na_{ij}\frac{\partial^2}{\partial x_i\partial x_j}-
\frac{\partial}{\partial t}\]
in $D$ to the function $u(0,x_2,...,x_n,t)$ we get
\begin{eqnarray}
\label{waldner}
\mathcal L\Big( u(0,x_2,...,x_n,t)\Big) &=&
\sum_{i,j=2}^na_{ij}(x_i,x_j,t)\frac{\partial^2}{\partial x_i\partial x_j}u(0,x_2,...,x_n,t)\\
\notag
&&-\frac{\partial}{\partial t}u(0,x_2,...,x_n,t)=0
\end{eqnarray}
since $u(0,x_2,...,x_n,t)$ satisfies the $(n-1)$-dimensional Black-Scholes equation
in the face $x_1=0$. Thus we have
\begin{equation}
\label{eq2}
\mathcal L\Big(u(x_1,x_2,...,x_n,t)-u(0,x_2,...,x_n,t)-\delta v(x_1,x_2,...,x_n,t)\Big)
\leq 0\mbox{ in $D$}.
\end{equation}
Applying the weak maximum principle to the inequalities (\ref{eq1}) and (\ref{eq2})
we find that
\[u(x_1,x_2,...,x_n,t)-u(0,x_2,...,x_n,t)-\delta v(x_1,x_2,...,x_n,t)\geq 0\mbox{ in $D$.}\]
Using $v(P^\prime)=0$ we get
\[u_{x_1}(P^\prime)\geq \delta v_{x_1}(P^\prime)=\delta>0\]
which finishes the proof.
\end{proof}

\begin{remark}
Note that the assumption that $\alpha_{ij}$ is a function merely of $x_i$ and $t$ (condition (ii) in
Hypothesis~\ref{hyp}) is
essential. Indeed, if  $\alpha$ instead would depend on the whole vector
$x$ and $t$, then the equality (\ref{waldner}) would not be true in general.
\end{remark}

\begin{remark}
Theorem~\ref{degenhopf} also holds if one includes lower order terms in the differential operator
$\mathcal L$. In addition to the assumptions on the coefficients specified in Theorem~\ref{degenhopf2}
below, one also needs to assume that $c$ and $b_i$, $i=2,...,n$ are independent of $x_1$, where
$c$ and $b_i$ are as defined in that theorem.

The theorem also generalizes to any non-tangential direction with non-positive time-component.
In that case, $N$ has to be chosen strictly larger than 1.
\end{remark}

To illustrate the result of Theorem~\ref{degenhopf} we give two examples.

\begin{example}{\bf (The Margrabe exchange option.)}
Let $n=2$, assume that
\[dX_i=\sigma_iX_i^{\beta_i/2}\,dW_i,\]
where $\sigma_i>0$ and $0\leq\beta_i<2$ for $i=1,2$, and let $g(x_1,x_2)=(x_2-x_1)^+$.
From Theorem~\ref{degenhopf} it follows that $u_{x_2}(x_1,0,t)>0$ for $x_2>0$ and $t>0$.
To investigate the derivative at the boundary $x_1=0$, note that adding the affine function $x_1$
to the contract function $g$ gives $(x_2-x_1)^++x_1=\max\{x_1,x_2\}$.
By Theorem~\ref{degenhopf} the value $\tilde u$
of this contract has a spatial derivative $\tilde u_{x_1}(0,x_2,t)>0$.
It follows that $u_{x_1}=\frac{\partial}{\partial x_1}(\tilde u-x_1)>-1$ for points $(0,x_2,t)$.

On the other hand, if $\beta_1=\beta_2=2$, i.e. if $X_1$ and $X_2$ are geometric Brownian motions,
then standard formulas for the value of the exchange option (compare \cite{M}) can be used to show
that $u_{x_2}(x_1,0,t)=0$ and $u_{x_1}(0,x_2,t)=-1$.
\end{example}

\begin{example}
{\bf (The call option in a CEV-model.)}
Let $n=1$, $g(x)=(x-K)^+$ and
\[dX=\sigma X^{\beta/2}\,dW\]
for some constants $K>0$, $\sigma>0$ and $\beta\in[0,2]$.
If $\beta=2$, then $X$ is a geometric Brownian motion and we know that $u_x(0,t)=0$ (compare
the example in the introduction), and if $\beta=1$, then one can use the valuation formula
in \cite{CR} (p. 161, equation (36)) to explicitly calculate
that $u_x(0,t)=\exp\{-\frac{2K}{\sigma^2t}\}>0$ for $t>0$.
Theorem~\ref{degenhopf} tells us that also in the remaining cases, i.e. for all $\beta<2$,
we have $u_x(0,t)>0$ for $t>0$.
\end{example}

In dimension $n=1$ we have the following result that for instance covers the example given
in the introduction. It shows that the condition that $\alpha_{11}\geq Cx_1^\beta$
for some $\beta\in(0,2)$ in Theorem~\ref{degenhopf} cannot be substantially weakened.

\begin{theorem}
\label{neghopf}
Let $n=1$, and assume that $a(x,t)\leq Cx^2$ for all $x$ and $t$ for
some constant $C$. Also assume that $g^\prime$ exists at $x=0$. Then $u_x(x,t)$ exists at $x=0$
and $u_x(0,t)=g^\prime(0)$ for all $t$.
\end{theorem}

\begin{proof}
By subtracting a constant and a suitable multiple of $x$ we may without loss of generality
assume that $g(0)=0$ and $g^\prime(0)=0$ (note that all affine functions $w$ satisfy
$\mathcal Lw=0$). Let $\ep>0$ be given.
Since $g$ is of at most polynomial growth, we can find $C_1>0$ and $N>1$ such that
\[g(x)\leq \ep x+C_1x^N\]
for all $x$. Note that $C_1$ depends on $\ep$, whereas $N$ can be held fixed if varying $\ep$.
Define the function $v$ by
\[v(x,t):=\ep x + C_1x^N + C_2tx^N\]
for some constant $C_2>0$ to be chosen. Then
\begin{eqnarray*}
v_t-a(x,t)v_{xx} &\geq& C_2x^N-Cx^2 C_1N(N-1)x^{N-2}-Cx^2C_2tN(N-1)x^{N-2}\\
&=& (C_2-C(C_1+C_2t)N(N-1))x^N,
\end{eqnarray*}
so for times $0\leq t\leq\frac{1}{2C(N^2-N)}=:t_0$ we have
\[v_t-a(x,t)v_{xx} \geq (C_2/2-CC_1N(N-1))x^N\geq 0\]
if $C_2$ is chosen large enough. Thus $v$ is a supersolution to the Black-Scholes equation
satisfying $v(x,0)\geq g(x)$. It follows from the maximum principle that $u\leq v$ for all
times $0\leq t\leq t_0$. Thus $u_x(0,t)\leq v_x(0,t)=\ep$ for such $t$. The same argument applied
to $-u$ gives $u_x(0,t)\geq -v_x(0,t)=-\ep$. Letting $\ep\to 0$ we arrive at
$u_x(0,t)=0$ for times $t\leq t_0$. Viewing
$u(x,t_0)$ as the initial condition the above argument shows that $u_x(0,t)=0$
also for $t_0\leq t\leq 2t_0$ (note that the same constant $N$ can be used again)
and thus by iteration also for all $t\geq 0$, which finishes the proof.
\end{proof}

There is no immediate generalization of Theorem~\ref{neghopf} to higher dimensions.
For instance, the price $u(x,t)$
in a geometric Brownian motion model of a call option on the sum of two assets
(i.e. $g(x_1,x_2)=(x_1+x_2-K)^+$) satisfies $u_{x_1}(0,x_2,t)\not= g_{x_1}(0,x_2)$.
On the other hand, the example with the exchange option shows
that the condition on $a_{11}$ in Theorem~\ref{degenhopf} is sharp also in
several dimensions.

Note that in the version of the Hopf boundary point lemma provided above (Theorem~\ref{degenhopf})
it is not assumed that $P^\prime$ is a minimum of the function $u$ but rather of the function
$u(x_1,x_2,...,x_n,t)-u(0,x_2,...,x_n,t)$. This is appropriate in the case when the boundary
conditions are defined inductively by solving the equation in lower dimensional faces.
Our proof also works in the situation where $u(x_1,x_2,...,x_n,t)$ has a minimum
in $P^\prime$. This situation is not the typical one when the boundary conditions are defined
as in the introduction, but it can of course occur in other types of degenerate
initial-boundary value problems.

\begin{theorem}
\label{degenhopf2}
Assume that $a$ is continuous and parabolic in $(0,\infty)^n\times (0,\infty)$, i.e.
$\xi a(x,t)\xi^*>0$ for all $\xi\in\R^n\setminus\{0\}$ and for all $x,t$.
Assume also that $u$ satisfies
\begin{equation}
\label{ineq}
\sum_{i,j=1}^n a_{ij}u_{x_ix_j}+\sum_{i=1}^nb_iu_{x_i}+cu-u_t\leq 0
\end{equation}
in $(0,\infty)^n\times (0,\infty)$, where $b_i=b_i(x,t)$ and $c=c(x,t)$ are continuous functions.
Let a point $P^\prime=(0,x^\prime_2,...,x^\prime_n,t_0)$
with $x^\prime_2,...,x^\prime_n,t_0>0$ be given, and assume that in a neighborhood of
$P^\prime$ we have
\[a_{11}\geq Cx_1^\beta,\]
\[a_{ii}\leq C\]
for $i=2,...,n$,
\begin{equation}
\label{b}
b_1\geq -C x_1^{\beta-1 +\delta},
\end{equation}
\[\vert b_i\vert\leq Cx_1^{\beta-2+\delta}\]
for $i=2,...,n$ and
\begin{equation}
\label{c}
c\geq -Cx_1^{\beta-2+\delta}
\end{equation}
for some constants $C>0$, $\delta>0$ and
$\beta\in[0,2)$. Also, in the same neighborhood, assume that $u(P)>u(P^\prime)$ for all points $P$
with $x_1>0$ and $t\leq t_0$. Then $u_{x_1}(P^\prime)>0$.
\end{theorem}

\begin{proof}
The proof follows along the lines of the proof of Theorem~\ref{degenhopf}
and is therefore omitted.
\end{proof}

As is shown in the next example, the lower bounds on $b$ and $c$ in Theorem~\ref{degenhopf2}
are sharp in the sense that the result fails for $\delta=0$.

\begin{example}
Let $n=1$. The function $u(x,t):=x^2/2$ is the solution of
\[\left\{\begin{array}{l}
u_t=x^\beta u_{xx}-x^{\beta-1}u_x\\
u(x,0)=x^2/2.\end{array}\right. \]
Note that (\ref{b}) is not satisfied and that $u_x(0,t)=0$.
Moreover, the same function $u$ also is the solution of
\[\left\{\begin{array}{l}
u_t=x^\beta u_{xx}-2x^{\beta-2}u\\
u(x,0)=x^2/2.\end{array}\right. \]
For this system (\ref{c}) is not satisfied.
\end{example}

\begin{remark}
If $\beta=0$ in Theorem~\ref{degenhopf2},
i.e. if the differential operator is uniformly parabolic, then one only needs the
inequality (\ref{ineq}) to hold in a parabolic frustrum, compare Lemma 2.6, p. 10 in \cite{L}.
Examining the proof of the Hopf lemma for non-uniformly operators, it is clear that this
also is true for $\beta\in[0,1)$, since $N$ in that case can be chosen to be 1.
\end{remark}

\section{The spatial derivative in dimension one}
\label{regularity}

In this section we use preservation of convexity for $n=1$ to deduce the existence of the derivative
$u_x(0,t)$ and also some regularity properties of this function as a function of time.
If the operator is uniformly parabolic, then $u_x(0,t)$ exists and is continuous for $t>0$
(even if $g^\prime(0)$ does not exist). To see this, note first that we can assume, without loss
of generality, that $g(0)=0$. Then, extending $g$ to a continuous odd function on $\R$ and
$\alpha$ to an even function in $x$, standard interior regularity results for parabolic PDE:s
yield that $u_x(0,t)$ is continuous. Note that, indeed, the solution to the extended problem
agrees with the solution to the original problem for non-negative $x$, since the extended solution
is odd in $x$ and thus vanishes at $x=0$.

We start this section with an example that shows that this is not true in
general for degenerate operators.

\begin{example}
{\bf (Power options in a geometric Brownian motion model.)}
Let $a(x,t)=x^2$ and $g(x)=x^\gamma$ for some constant $\gamma>0$.
Then it is easy to check that
\[u(x,t)=x^\gamma \exp\{(\gamma^2-\gamma)t\}.\]
Thus $u_x(0,t)$ does not exist if $0<\gamma<1$.
\end{example}

\begin{theorem}
\label{usc}
Assume that the contract function $g$ is convex and that $g^\prime(0+)>-\infty$.
Then the derivative $u_x(0,t)$ exists for all $t$.
Moreover, the function $t\mapsto u_x(0,t)$ is increasing
and upper semi-continuous.
\end{theorem}

\begin{proof}
Recall, see \cite{JT1}, that convexity is preserved and that the option price increases
in time to maturity, i.e. $x\mapsto u(x,t)$ is convex for each fixed $t\geq 0$ and
$t\mapsto u(x,t)$ is increasing for each fixed $x$, see also \cite{BGW}, \cite{EKJPS} and \cite{H}.
Since $u(0,t)=g(0)$ it follows that $u_x(0,t)$ exists and that $u_x(0,t)$ is increasing.
Moreover, using the continuity in $t$ and the spatial convexity of $u$,
it follows that $u_x(0,t)$ is upper semi-continuous.
\end{proof}

\begin{remark}
Note that $u_x(0,t)$ also exists if the contract function can be written
as a difference of two convex functions (both with finite derivative in the origin).
Also note that the example above with $g(x)=x^\gamma$ where $0<\gamma<1$ (or rather $g(x)=-x^\gamma$)
shows that the assumption in Theorem~\ref{usc} about $g^\prime(0)>-\infty$ is essential.
\end{remark}

\begin{remark}
In higher dimensions convexity is in general no longer preserved, compare \cite{EJT}.
Therefore the above proof of Theorem~\ref{usc} is not applicable to problems with several
underlying assets. It remains an open question to determine conditions under which
the spatial derivatives exist at the boundary if $n\geq 2$.
\end{remark}

The conclusions of Theorem~\ref{usc} that can be drawn more or less immediately from
the preservation of convexity cannot be improved in the sense that $u_x(0,t)$ need not
be continuous in $t$. Indeed, we end this article with the construction of an example
where $\alpha$ is
continuous and locally Lipschitz in $x$, the contract function $g$ is convex and in
$C^\infty([0,\infty))$ and yet $u_x(0,t)$ fails to be lower semi-continuous.

\begin{example} {\bf ($u_x(0,t)$ need not be continuous as a function of $t$).}
The example is constructed by patching together two functions
$v:[0,\infty)\times[0,t_0)\to [0,\infty)$ and $w:[0,\infty)\times[t_0, \infty)\to [0,\infty)$
for some given $t_0>0$.

To define $v$, let $h\in C^\infty([0,\infty))$ be convex, non-negative and satisfy
$h(0)=0$, $h^\prime(0)=0$, and $h^\prime(y)=1$ and $h^{\prime\prime}(y)=0$ for
$y\geq y_0>0$. Let
\[C:=\lim_{y\to\infty}(yh^\prime(y)-h(y))=y_0-h(y_0).\]
Define $v:[0,\infty)\times [0,t_0)\to [0,\infty)$ by
\[v(x,t)=x^2e^t+(t_0-t)h(\frac{x}{t_0-t}).\]
It follows that
\[v_t=\tilde a(x,t)v_{xx}\]
in $(0,\infty)\times(0,t_0)$ where
\[\tilde a(x,t)=\frac{x^2e^t(t_0-t)+xh^\prime(\frac{x}{t_0-t})-(t_0-t)h(\frac{x}{t_0-t})}
{2e^t(t_0-t)+h^{\prime\prime}(\frac{x}{t_0-t})}.\]
Note that $v_x(0,t)=0$ for $0\leq t<t_0$, and for $x>0$ we have
\[\lim_{(y,t)\to(x,t_0)}v(x,t)=x^2e^{t_0}+x,\]
\[\lim_{(y,t)\to(x,t_0)}v_t(x,t)=x^2e^{t_0}+C,\]
\[\lim_{(y,t)\to(x,t_0)}v_x(x,t)=2xe^{t_0}+1\]
and
\[\lim_{(y,t)\to(x,t_0)}v_{xx}(x,t)=2e^{t_0}.\]
Next define $w$ as the unique solution to
\[\left\{\begin{array}{ll}
w_t=\frac{x^2+Ce^{-t_0}}{2}w_{xx} & \mbox{for $(x,t)\in(0,\infty)\times(t_0,\infty)$}\\
w=0 & \mbox{if $x=0$}\\
w=x^2e^{t_0}+x & \mbox{if $t=t_0$.}\end{array}\right.\]
Then it is straightforward to check that
\[w_t(x,t_0)=x^2e^{t_0}+C,\]
\[w_x(x,t_0)=2xe^{t_0}+1\]
and
\[w_{xx}(x,t_0)=2e^{t_0}\]
for $x>0$.
It follows that
\[u(x,t)=\left\{\begin{array}{ll}
v(x,t) & \mbox{if $0\leq t<t_0$}\\
w(x,t) & \mbox{if $t_0\leq t<\infty$}\end{array}\right.\]
solves
\[u_t=a(x,t)u_{xx}\]
where
\[a(x,t):=\left\{\begin{array}{ll}
\tilde a(x,t) & \mbox{if $0\leq t<t_0$}\\
\frac{x^2+C}{2} & \mbox{if $t_0\leq t<\infty$.}\end{array}\right.\]
Note that $\alpha:=\sqrt{2a}$ is continuous on $(0,\infty)\times(0,\infty)$, locally Lipschitz in
the $x$-variable and that $\alpha$ satisfies the growth condition
\[\alpha(x,t)\leq D(1+x)\]
for some positive constant $D$.
Recall that option prices with convex contract functions are increasing
in the time to maturity, compare \cite{JT1}, i.e. the function $t\mapsto w(x,t)$ is increasing.
Thus $w_x(0,t)\geq 1$ for $t\geq t_0$. Therefore the function $u_x(0,t)$ is not continuous.
\end{example}

\end{document}